\documentclass[11pt]{amsart}
\usepackage{geometry}                % See geometry.pdf to learn the layout options. There are lots.
\usepackage{graphicx}
\usepackage{amssymb}
\usepackage{epstopdf}
\usepackage{hyperref}
\usepackage{mathrsfs}
\usepackage[arrow,curve,matrix]{xy}
\input xy
\newcommand{\Gr}{{\mathrm{Gr}}}

\newcommand{\rk}{{\rm rk}}

\newcommand{\im}{{\mathrm{im}}}
\newcommand{\pr}{{\mathrm{pr}}}

\theoremstyle{plain}

\newtheorem{thm}{Theorem}[section]
\newtheorem{prop}[thm]{Proposition}
\newtheorem{lem}[thm]{Lemma}

\theoremstyle{definition}
\newtheorem{defn}[thm]{Definition}
\newtheorem{ex}[thm]{Example}
\newtheorem{rmk}[thm]{Remark}

\def\CC{{\mathbb C}}
\def\FF{{\mathbb F}}
\def\GG{{\mathbb G}}
\def\PP{{\mathbb P}}

\def\RR{{\mathbb R}}

\def\cO{\mathcal{O}}

\def\cO{\mathcal{O}}

\def\cQ{\mathcal{Q}}

\def\cS{\mathcal{S}}

\def\fQ{\mathfrak{Q}}
\def\fS{\mathfrak{S}}

\title[Holomorphic vector bundles and totally geodesic foliations]{Holomorphic vector bundles on K\"ahler manifolds and totally geodesic foliations on Euclidean open domains}

\author{Monica Alice Aprodu}
\address{"Dun\u area de Jos" University, Str. Domneasc\u a, nr. 47, 800008 Gala\c ti, Romania}
\address{"Simion Stoilow" Institute of Mathematics of the Romanian Academy, P.O. Box 1-764, 014700 Bucharest, Romania}
\email{Monica.Aprodu@ugal.ro}
\urladdr{www.math.ugal.ro}

\author{Marian Aprodu}
\address{"Simion Stoilow" Institute of Mathematics of the Romanian Academy, P.O. Box 1-764, 014700 Bucharest, Romania}
\email{Marian.Aprodu@imar.ro}
\urladdr{www.imar.ro/$\sim$aprodu}

\thanks{Monica Alice Aprodu was partly supported by the CNCS grant PNII-ID-PCE-2011-3-0362. Marian Aprodu was partly supported by the Alexander von Humboldt Foundation and Max Planck Institut f\"ur Mathematik Bonn; the author is grateful to Humboldt Universit\"at zu Berlin and Max Planck Institut f\"ur Mathematik Bonn for hospitality.}
\date{\today}                                           % Activate to display a given date or no date
\begin{document}
\maketitle

\begin{abstract}
In this Note we establish a relation between sections in globally generated holomorphic vector bundles on K\"ahler manifolds, isotropic with respect to a non-degenerate quadratic form, and totally geodesic foliations on Euclidean open domains. We find a geometric condition for a totally geodesic foliation to originate in a holomorphic vector bundle. For codimension-two foliations, this description recovers \cite{BW} Theorem 2.8 and Theorem 2.15. The universal objects that play a key role are the orthogonal Grassmannians.
\end{abstract}

%\section{Introduction}

In \cite{BW}, P. Baird and J. C. Wood have solved completely a problem posed by Jacobi \cite{Jac}. In modern language, the question was to classify harmonic morphisms with one-dimensional fibres from open domains in $\RR^3$. In other words, it was asked to describe (locally) foliations in lines on open domains in $\RR^3$ which verify a geometric extra--condition, horizontal weak conformality (HWC). Baird and Wood interpreted horizontal weak conformality as a holomorphic variation of the leaves. More precisely, assuming the foliation was simple, they associate a meromorphic map from the leaf space to a complex projective quadric and indicate an inverse construction. Two distinct cases are found, corresponding to whether or not the leaves pass through the origin.Their approach works for codimension--two foliations on higher--dimensional Euclidean open domains.

In this Note, we extend the results of Baird and Wood to foliations of higher codimension. We establish a relation between holomorphic maps to an orthogonal Grassmannian, together with a section in the universal bundle, and certain foliations with totally geodesic leaves on Euclidean open domains, Theorem \ref{thm:holo-foli}, Theorem \ref{thm:foli-holo}. A key feature of the orthogonal Grassmannian is the universal property, it represents the functor of isomorphism classes of isotropic quotient bundles of trivial bundles, Proposition \ref{prop:univ} in section \ref{sec:OG}. The definition and the basic properties of the orthogonal Grassmannians are recalled also in section \ref{sec:OG}.

\medskip

\noindent{\em Notation and convention.} All the manifolds are considered without boundary. If $M$ is a real manifold, $T_\RR M$ denotes the real tangent bundle, and $T_\CC M:=T_\RR M\otimes_\RR\CC$ denotes the complexified tangent bundle. For a complex manifold $N$, we denote by $T^+N$ the holomorphic tangent bundle. Recall that, as {\em real} vector bundles, $T_\RR N$ and $T^+N$ are isomorphic, although they differ as subbundles in $T_\CC N$. If $\varphi:M\to N$ is a submersion between Riemann manifolds and $X$ is a vector field on $M$, we denote by $\mathscr{H}X$ and $\mathscr{V}X$ the horizontal, respectively the vertical components of $X$. If $E$ is a vector bundle on a manifold $M$, we denote by $\Gamma(E)$ the space of $\mathscr{C}^\infty$--sections in $E$.

\section{Orthogonal Grassmannians and isotropic quotient bundles}
\label{sec:OG}

The general theory of orthogonal Grassmannian is developed for subspaces of a given vector space \cite[p. 387]{FH}, \cite{Raw} etc. In this section, we adapt the theory to work with quotients, and the reason is twofold. Holomorphic subbundles of trivial bundles rarely admit holomorphic sections, whereas quotient bundles have plenty of sections. From the sheaf-theoretic viewpoint, from a quotient locally free sheaf we can always recover a subsheaf by dualization, whereas dualizing a locally free subsheaf does not necessarily define a locally free quotient. With this view-point, the normal bundle of a foliation is a quotient bundle, rather than a subbundle.%\marginpar{\tiny de reformulat}

\subsection{Isotropic quotient spaces} 
Let $V$ be a $m$--dimensional real vector space endowed with an inner-product $h$. Denote by $V^*$ real dual space, and $V_\CC^*$ the complex-dual. Put $h_\CC$ the symmetric bilinear form induced on the complexification $V_\CC=V\otimes_\RR\CC$. The musical isomorphism $\flat:V\to V^*$, $v\mapsto v^\flat:=h(v,-)$ induces an inner-product $h^*$ on $V^*$ and a symmetric bilinear form $h^*_\CC$ on~$V_\CC^*$.

\begin{defn}
An $n$--dimensional (complex) quotient space $V_\CC\to\Lambda$ is called {\em isotropic} if $h^*_\CC|_{\Lambda^*\times\Lambda^*}\equiv 0$, where the inclusion $\Lambda^*\subset V_\CC^*$ is obtained by duality.
\end{defn}

Isotropic quotients are $\RR$--quotients of $V$, too \cite{AA}:

\begin{lem}
\label{lem:isotropic}
Let $\Lambda$ be an isotropic quotient of $V_\CC$. The induced $\RR$--linear map $V\to \Lambda$ remains surjective. 
\end{lem}

\proof
The conclusion is proved by direct computation, using an orthonormal basis, see \cite[Lemma 3.2]{AA} for details.
\endproof

As a real vector space, an isotropic quotient $\Lambda$ carries a complex structure, compatible with the inner-product induced by $h^*$. Conversely, any $2n$--dimensional quotient of $V$, together with a complex structure compatible with $h^*$ corresponds to an isotropic quotient of $V_\CC$. By abusing terminology, we shall sometimes refer to $\Lambda$ as an {\em isotropic quotient} without mentioning whether it is regarded as a quotient of $V$ or of $V_\CC$.

\begin{rmk}
The following alternate interpretation is particularly useful for our purposes. It was observed in \cite{Raw} that there is a corespondence between isotropic quotients and skew-symmetric endomorphisms $F$ of $V$ of rank $2n$ with the condition $F^3+F=0$, given as follows. To any endomorphism $F$, one associates the quotient space $\im(F)=V/\ker(F)$, together with the induced complex structure. Skew-symmetry ensures that this complex structure is compatible with the inner-product. Conversely, to any $2n$--dimensional quotient space $\varphi:V\to \Lambda$ with a compatible complex structure $J$, one associates naturally a skew-symmetric endomorphism $F$ of $V$ of rank $2n$. The inner--product $h$ induces an orthogonal decomposition $V=\ker(\varphi)\oplus \Lambda$, and $F$ is defined to be identically zero on $\ker(\varphi)$ and to be $J$ on~$\Lambda$. 
\end{rmk}

\subsection{The orthogonal Grassmannian and universal bundles}
For any $1\le n\le m/2$,  the {\em orthogonal Grassmannian} $\FF=F_n(V^*,h^*)$, is defined as the parameter space of $n$-dimensional isotropic quotient spaces of $V_\CC$, i.e.
\[
F_n(V^*,h^*):=\{\Lambda^* \subset V^*_\CC:\ h^*_\CC|_{\Lambda^*\times\Lambda^*}\equiv 0\}
\]

As a complex subvariety of the Grassmannian $\GG=\Gr_n\left(V^*_\CC\right)$, the variety $F_n(V^*,h^*)$ it is defined by quadratic equations, whence the alternate name, {\em Grassmann quadric}~\cite{AA}, \cite[p. 97]{Raw}. In other references, the variety $\mathbb F$ is denoted by $OG(n,m)$. It identifies with the Fano scheme of $(n-1)$--dimensional projective subspaces of $\PP(V^*_\CC)$ contained in the quadric defined by the quadratic form associated to $h^*_\CC$~\cite{Harris}. It is interesting to note that Fano schemes are one of the few classes of Hilbert schemes which are smooth.

\medskip

The orthogonal Grassmannian comes with a universal short exact sequence of vector bundles: 
\begin{equation}
\label{eqn:complex}
0\to \fS \to \FF\times V_\CC\to \fQ\to 0,
\end{equation}
where $\fQ$ is the restriction to $\FF$ of the universal rank--$n$ quotient bundle from the Grassmannian~$\GG$. Let $\cQ=\, _\mathbb{R}\fQ$ be the {\em realification} of $\fQ$, i.e. the bundle $\fQ$ seen as real vector bundle. Applying Lemma \ref{lem:isotropic}, we obtain a short exact sequence of real vector bundles
\begin{equation}
\label{eqn:real}
0\to \cS \to \FF\times V\to \cQ\to 0,
\end{equation}
where $\cS:=\fS\cap (\FF\times V)$ is a real bundle of rank $(m-2n)$. 
$\cQ$ is the prototype of an {\em isotropic} quotient bundle, in the following sense.

\begin{defn}
\label{defn:isotropic quot}
Let $M$ be a real manifold and $V$ be a $m$--dimensional real vector space with an inner-product $h$. A quotient bundle $M\times V\to E$ of rank $2n$ is called {\em isotropic} if it admits an almost complex structure, i.e. an endomorphism $J_E$ with $J_E^2=-\mathrm{id}$, such that for any $x\in M$ the fibre $E_x$ is an isotropic quotient of $V$. If $M$ is a complex manifold and $J_E$ is compatible with the complex structure of $M$, then $E$ is called a {\em isotropic holomorphic quotient}.
\end{defn}

The compatibiliy relation for isotropic holomorphic quotients has the following meaning. Choose {\em local} holomorphic trivializations $E\cong M\times \mathbb C^n$. Then the quotient map $M\times V\to E$ corresponds to a smooth map $M\times V\to \mathbb C^n$ obtained by projection onto $\mathbb C^n$.  This map has to be holomorphic on $M$. Since the $\mathbb R$--linear maps $V\to \mathbb C^n$ extend to $\mathbb C$--linear maps $V_{\mathbb C}\to \mathbb C^n$, it follows that isotropic holomorphic quotients are realifications of holomorphic quotients of trivial holomorphic bundles.

\begin{rmk}
\label{rmk:isotropic quot}
Isotropic quotients can be defined also for arbitrary real vector bundles endowed with a smooth metric, for instance tangent bundles of Riemann manifolds. 
\end{rmk}

\subsection{The universal property}
The next result follows directly from the definition and uses the universal property of the Grassmannian $\Gr_n(V_\CC)$. It settles the place of the orthogonal Grassmannian in the category of manifolds and explains our preference for quotient spaces over subspaces.

\begin{prop}
\label{prop:univ}
The orthogonal Grassmannian represents the contravariant functor $(\underline{\mathrm{Mfld}})\to (\underline{\mathrm{Sets}})^{op}$ from the category of real manifolds to the category of sets that associates to any manifold $M$ the set of isomorphism classes of almost--complex isotropic quotient bundles of rank $2n$ of $M\times V$. Restricted to the subcategory of complex manifolds, it represents the functor that associates to any complex manifold $N$, the set of isomorphism classes of isotropic holomorphic bundle quotients of $N\times V_\CC$. In other words, any almost-complex isotropic quotient bundle is the pull-back of $\cQ$ through a smooth map, and any isotropic holomorphic quotient is the pull-back of $\fQ$ through a holomorphic map.
\end{prop}

Since the smooth quadric in the $(m-1)$--dimensional complex projective space $\PP(V_\CC)$ coincides with $F_1(V^*,h^*)$, Proposition \ref{prop:univ} implies that the quadric parametrizes isotropic holomorphic rank-one quotients, which is clear.

\subsection{The orthogonal Grassmannian as a homogeneous manifold}
\label{sec:homog}
Fixing $F\in \mathbb F$ and viewing it as a skew-symmetric endomorphism of $V$, one obtains the description of $\mathbb F$ as the quotient of the orthogonal group $\mathscr{O}(V,h)$ through the stability group $\{g\in \mathscr{O}(V,h):\ gFg^{-1}=F\}$. This is proved interpreting $\FF$ as a subset of the Lie algebra 
\[
o(V,h) = \{ A\in\mathrm{End}(V)\ :\ A+A^\top =0\}
\] 
and observing that it forms a single conjugacy class under the adjoint action of $\mathscr{O}(V,h)$, see \cite[p. 97]{Raw}. This description has the advantage of giving rise to a simple description of the holomorphic tangent space $T^+_F\mathbb F$ in terms of matrices, as follows. If $\mathrm{ad}(F):o(V,h)\to o(V,h)$ denotes the bracketing with $F$, i.e. $\mathrm{ad}(F)(A)=[A,F]=AF-FA$ for any $a\in o(V,h)$, then $\im(\mathrm{ad}(F))=[F,o(V,h)]$ identifies with $T_F\mathbb F$ inside $o(V,h)$ \cite[Proposition 3.1]{Raw}. A direct computation shows that $\mathrm{ad}(F)|_{T_F\mathbb F_{\mathbb C}}$ satisfies the equation 
\[
\mathrm{ad}(F)^4+5\mathrm{ad}(F)^2+I=0,
\] 
which shows that its eigenvalues can only be $\pm2i$ and $\pm i$. Then $T^+_F\mathbb F$ is the direct sum of the eigenspaces corresponding to the eigenvalues $+2i$ and $+i$, and a multilinear algebra argument  implies that 
\begin{equation}
\label{eqn:T+F}
T^+_F\mathbb F=\{ A\in o(V,h)_{\mathbb C}\ :\ AV^+=0,\ AV^0\subset V^+,\ AV^-\subset V^0+V^+ \},
\end{equation}
where $V^+$, $V^-$ and $V^0$ are the eigenspaces of $V_{\mathbb C}$ corresponding to the eigenvalues $+i$, $-i$ and, respectively, $0$ of $F$ \cite[Proposition 3.3]{Raw}.

There is a natural submersion $\Phi$ from $\mathbb F$ the Grassmannian $\widetilde{\mathrm{Gr}}_{2n}(V)$ of oriented $2n$--dimensional real quotients of~$V$, which forgets the complex structure on $\im(F)$. The differential $\mathrm{d}\Phi$ has the following simple description. To any $F\in\mathbb F$ and any $B=FA-AF\in T_F\mathbb F$ with $A\in o(V,h)$ we associate the map $(FA)|_{\ker(F)}:\ker(F)\to \im(F)$. Note that  this map coincides with $B|_{\ker(F)}$, as $(AF)|_{\ker(F)}=0$. In the case $n=1$, the differential $\mathrm{d}\Phi$ is an isomorphism, and the induced complex structure on the space of $\mathbb R$--linear maps from $\ker(F)$ to $\im(F)$ is given by the composition with $F$; indeed in this case, one verifies directly that $\mathrm{ad}(F)^2+I=0$. In general, for $n$ arbitrary, we have the following useful identity
\begin{equation}
\label{eqn:dPhi}
\mathrm{d}\Phi(J_{\mathbb F}(B))=F\circ \mathrm{d}\Phi(B),
\end{equation}
for any $B=FA-AF\in T_F\mathbb F$, where $J_{\mathbb F}$ denotes the complex structure on $\mathbb F$.

\section{$f$--structures}
\label{sec:F}

\subsection{$f$--structures and isotropic quotients}
The isotropy condition is effectively used in the theory of harmonic maps \cite{Lou}, \cite{BW-book}, \cite{AAB}, \cite{AA}. The easiest case is the following:

\begin{lem}
\label{lem:PHWC}
Let $V$ be an $m$--dimensional real vector space with an inner-product $h$, $\Lambda$ be a complex vector space and $\tau: V\to \Lambda$ be a surjective $\RR$--linear map. Then $\Lambda$ is an isotropic quotient if and only if it pulls-back germs of holomorphic functions on $\Lambda$ to germs of harmonic functions on $V$.
\end{lem}

\proof
Note that an $\mathbb R$--linear map is harmonic. To conclude, apply \cite[Proposition 2]{Lou}.
\endproof

The latter condition is usually called in literature {\em PHWC} or {\em pseudo-horizontally weakly conformal} \cite{Lou}, \cite{BW-book}, the reason being that it was introduced as an extension of horizontally weakly conformal maps. 

An immediate consequence of Lemma \ref{lem:PHWC} and of \cite{Lou} is the fact  that complex quotients of complex vector spaces are isotropic. (this can also be proved directly.)

Note that the fibres of a linear map as in Lemma \ref{lem:PHWC} are affine subspaces. It is the first and the simplest sample of a foliation on an Euclidean domain that originates in a holomorphic vector bundle. The bundle in question is the trivial bundle on $\Lambda$ with fibre $\Lambda$ and is the pullback of the universal bundle $\fQ$ from $\Lambda$ to $F_n(V^*,h^*)$ through a constant map.

We recall the following \cite{IY}, see also \cite{Raw}:

\begin{defn}[Yano]
\label{defn:F}
Let $(M^m,g_M)$ be a Riemann manifold. An $f$-structure of rank $2n\le m$ on $M$ is a skew-symmetric endomorphism $F$ of the tangent bundle $TM$ with $F^3+F=0$ and $\rk(\ker(F))=m-2n$.
\end{defn}

Skew-symmetry is required for compatibility with the Riemann metric. For $2n=m$, an $f$--structure is an almost complex structure compatible with the Riemann metric. As pointed out in \cite{Raw}, parallel $f$--structures, i.e. for which $\nabla F=0$, are subject to several restrictions.

If $F$ is any $f$--structure on $M$, then the bundle $TM\to \im(F)=TM/\ker(F)$ is isotropic in the sense of Remark \ref{rmk:isotropic quot}. Conversely, given an isotropic quotient of the tangent bundle of $M$, it defines an $f$--structure on $M$. Hence we have a one-to-one correspondence between $f$--structures and isotropic quotients of the tangent bundle. If the manifold is parallelizable, then they are quotients in the usual sense Definition \ref{defn:isotropic quot}.

$f$--structures induced by maps to Hermitian manifolds are of particular interest for us. Suppose that $(N^{2n},J_N,g_N)$ is a Hermitian manifold and $\varphi:M\to N$ is a submersive surjection. The differential the pullback of the real tangent bundle of $N$, $\mathrm{d}\varphi$ exhibits $\varphi^{-1}(T_\RR N)$ as a quotient of $T_\RR M$. 

\begin{prop}
\label{prop:PHWC}
The bundle $\varphi^{-1}(T_\RR N)$  is an isotropic quotient of $T_\RR M$ if and only if the map is pseudo horizontally weakly conformal.
\end{prop}

Hence, a PHWC submersion naturally yields to an $f$--structure on $M$, denoted by $F^\varphi$, see \cite[section 3]{Lou}, \cite[section 4]{LM}. For a precise definition of pseudo horizontal weak conformality, we refer to \cite[Definition 1]{Lou},~\cite{BW-book} and Proposition \ref{prop:PHWC} can be considered as an alternate definition.

\medskip

In the sequel, we will restrict to a special class of PHWC maps, called {\em sPHH} or {\em strongly pseudo-horizontally homothetic} \cite{AAB}, \cite{AA}. They are characterized by a partial parallel condition on the associated $f$--structure, namely $(\mathscr{H}\nabla F^\varphi)(X,Y)=0$ (equivalently $\mathscr{H}\nabla_XF^\varphi(Y)= F^\varphi(\nabla _XY)$) for any vector field $X$ and any horizontal vector field $Y$. If $\mathscr{H}\nabla F^\varphi$ vanishes on horizontal vector fields, then the map is called {\em pseudo-horizontally homothetic},  \cite{AAB}, \cite{AA}. We can easily generalize the sPHH (or the PHH) condition to any $f$--structure $F$ on $M$.

%This terminology points out to the following definition.
%
%\begin{defn}
%Let $(M^m,g_M)$ be a Riemannian manifold. An $f$--structure $F$ on $M$ is called {\em PHH}, respectively {\em sPHH}, if $(\nabla F)(X,Y)=0$ for any vector fields $X,Y\in\ker(F)^\perp$, respectively any vector field $X$ and any vector field $Y\in \ker(F)^\perp$.
%\end{defn}

\subsection{The Nijenhuis tensor}
The integrability of an $f$--structure in controlled by an associated tensor \cite{IY}. The $f$--structure $F$ on a Riemannian manifold $(M^m,g_M)$ is called {\em partially integrable} if $\ker(F)^\perp$ is integrable and the almost complex structure induced on each integrable manifold is integrable (i.e. they are complex manifolds) \cite[p. 220]{IY}.
Note that $\ker(F)=\im(F^2+I)$ and $\ker(F)^\perp=\ker(F^2+I)=\im(F)$. In analogy to almost complex structures, we have a Nijenhuis tensor:
\[
N_F(X,Y):=[FX,FY]-F[FX,Y]-F[X,FY]+F^2[X,Y].
\]
In \cite{IY} the authors prove that partial integrability is equivalent to the vanishing of $N_F$ on vector fields in $\im(F)=\ker(F)^\perp$.

\section{sPHH maps and sections in holomorphic bundles}

In this section, we obtain a local description of totally geodesic foliations induced by sPHH maps on Euclidean domains in terms of  sections in holomorphic vector bundles on the target manifolds. Similar techniques can be applied for foliations on other parallelizable Riemann manifolds, such a flat tori.

\subsection{sPHH maps from sections in holomorphic bundles}
Let $V$ be a $m$--dimensional real vector space with an inner-product $h$ and $1\le n\le m/2$ be an integer. 
Consider $N$ a complex manifold of complex dimension $n$ and $\xi :N\to F_n(V^*,h^*)$ a holomorphic map. The pullback of the universal sequence provides us with a universal quotient $\tau:N\times V\to \cQ$. Consider $\sigma$ a holomorphic section (possibly zero) of $\fQ$, and put $V_{\sigma}:=\tau^{-1}(\sigma(N))$. Note that $V_\sigma$ is an $m$--dimensional manifold and the first projection $\pr_1:V_\sigma\to N$ is an affine bundle. Identifying $N$ with $\sigma(N)$, this projection identifies with the restriction of $\tau$.

The  fibres of $\pr_1$ are mapped to $(m-2n)$--dimensional affine subspace of~$V$ through the second projection $\pr_2:V_\sigma\to V$. Choose $U\subset V$ an open set over which $\pr_2$ is a diffeomorphism. For any $x\in U$, there exists a single $(m-2n)$--plane through $x$ which is the image of the fibre $\pr^{-1}(y)$ for some $y\in N$. We obtain a well-defined submersion $\varphi :U\to N$ by mapping $x$ to $y$, such that $V_\sigma$ identifies (locally) with the graph of $\varphi$. For simplicity, let us assume that $\varphi$ is surjective.
To this setup, we associate an $f$--structure on $U$ seen as a manifold, section \ref{sec:F}.

If we suppose that $N$ is K\"ahler, then we can prove that the map $\varphi$ is sPHH.
Indeed, since the property sPHH is local on the base, we may assume that $\mathcal Q$ is trivial, isomorphic to $N\times \mathbb C^n$ and the section $\sigma$ corresponds to a holomorphic map $(\sigma_1,\ldots,\sigma_n):N\to \mathbb C^n$ and the morphism $\tau$ corresponds to a map $\tau:N\times V\to \mathbb C^n$, holomorphic on $N$ and $\mathbb R$--linear on $V$. By our assumptions, $\tau$ is moreover PHWC on $V$. Up to a permutation of factors, we can assume that $\sigma_1,\ldots,\sigma_i$ are identically zero, and $\sigma_{i+1},\ldots,\sigma_n$ are not identically zero for some $i$.
Note that if the map $\varphi$ is sPHH over a dense open subset of $N$ then it is sPHH everywhere, and hence we may assume furthermore that $\sigma_{i+1},\ldots,\sigma_n$ are nowhere vanishing.
Define the map $G:N\times V\to \mathbb C^n$ by 
\[
G(x,v)=\left(\tau_1(x,v),\ldots,\tau_i(x,v),\frac{1}{\sigma_{i+1}(x)}\tau_{i+1}(x,v),\ldots,\frac{1}{\sigma_{n}(x)}\tau_n(x,v)\right).
\]
The map $\varphi$ is a solution for the implicit equation
\[
G(\varphi(v),v)=(0,\ldots,0,1\ldots,1).
\]
Since $\tau$ is holomorphic on $N$ and linear and PHWC on $V$, so is $G$ and hence the hypotheses of \cite[Theorem 2.1, Remark 2.2]{AA} are verified, which implies that $\varphi$ is~sPHH.
\endproof

We have obtained the following (compare to \cite{BW}, Theorem 2.15, for the case $n=1$)

\begin{thm}
\label{thm:holo-foli}
To any holomorphic map from a given K\"ahler manifold $N$ to the orthogonal Grassmannian and a holomorphic section in the pullback of the universal quotient bundle one can associate a totally geodesic foliation on an Euclidean domain $U$. This foliation corresponds to a sPHH $f$--structure $F$ on $U$. \end{thm}

\begin{ex}
\label{ex:Pn}
Let $\PP=\PP^n_\CC$ be the $n$--dimensional complex projective space and consider the holomorphic vector bundle $\mathcal Q=T^+\PP\otimes _{\cO_{\PP}}\cO_{\PP}(-1)$, where $\cO_{\PP}(-1)$ is the universal line subbundle. It is well-known that $\mathcal Q$ is globally generated, a simple computation shows that it is actually an isotropic quotient of $\PP\times V$, where $V=H^0(\mathcal O_\PP(1))^*\cong \mathbb R^{2n}$. If $\sigma$ denotes the zero--section in $\mathcal Q$, then, with the previous notation, $V_\sigma\subset \PP\times V$ is the incidence variety, i.e. the total space of the holomorphic line bundle $\mathcal O_\PP(-1)$. The projection to $\PP$ realizes $V_\sigma$ as a vector bundle and the projection to $V$ is the blowup of the origin. 
The induced submersion is the natural map $\pi:V\setminus\{0\}\to \PP=\PP^n_\CC.$ Note that $\pi$ is not only sPHH, but a Riemannian submersion.
\end{ex}

In the next subsection, we show a converse of Theorem \ref{thm:holo-foli}.

\subsection{Sections in holomorphic bundles from foliations}
Let $U$ be a domain in an $m$--dimensional Euclidean space $V$ with an inner-product $h$ and denote, as before, $o(V,h)$ the Lie algebra of the orthogonal group $\mathscr{O}(V,h)$. Let  $\varphi$ be a sPHH submersion with totally geodesic fibres from $U$ to a K\"ahler manifold $N$ of complex dimension $n$. Denote by $F=F^\varphi$ the $f$--structure on $U$ induced by $\varphi$, section~\ref{sec:F}.

\begin{lem}
\label{lem:vertical}
For any vertical vector field $X$, we have $\nabla_XF=0$.
\end{lem}

\proof
If $Y$ is a horizontal vector field, then $(\nabla_XF)(Y)=\mathscr{V}\nabla_XFY$, from the sPHH assumption. For any vertical $Z$ we have
\[
g_M(\nabla_XFY,Z)+g_M(FY,\nabla_XZ)=Xg_M(FY,Z)=0
\]
and, since the fibres are totally geodesic, $\nabla_XZ$ is vertical implying $g_M(\nabla_XFY,Z)=0$. It follows that $(\nabla_XF)(Y)=0$.

If $Y$ is vertical then $\nabla_XY$ is vertical, since the fibres are totally geodesic, hence $F\nabla_XY=\nabla_XFY=0$ implying $(\nabla_XF)(Y)=0$. 
\endproof

We consider the smooth map $\Omega$ from $U$ to the orthogonal Grassmannian $\mathbb F$ that associates to any $x\in U$ the horizontal space at $x$. In other words, seeing $\mathbb F\subset o(V,h)$, $\Omega$ maps a point $x\in U$ to $F_x:V\cong T_xU\to V\cong T_xU$.

Lemma \ref{lem:vertical} implies that the map $\Omega$ is constant along the fibres of $\varphi$ i.e. it factors through a map $\omega:N\to \mathbb F$. 

Using the universal property, it follows that $\ker(F)$ is a pull-back of a bundle on $N$. Under these assumptions, the natural short exact sequence 
\begin{equation}
\label{eqn:exact}
0\to \ker(F)\to T_\RR U\to \varphi^{-1}(T_\RR N)\to 0
\end{equation}
descents to a short exact sequence on $N$ such that that the quotient $Q^\varphi=\omega^*(\mathcal Q)$ is isotropic. Note that $\varphi^{-1}(Q^\varphi)$ coincides with $\varphi^{-1}(T_\RR N)$. 

Denote by $\tau$ the natural map $N\times V\to Q^\varphi$ and consider the graph $\Gamma_\varphi\subset N\times V$ of $\varphi$. The image of $\Gamma_\varphi$ via $\tau$ is a submanifold of $Q^\varphi$. It is the image of a section $\sigma$ of $Q^\varphi$ which has the following interpretation. Over each $y\in N$, the fibre $\varphi^{-1}(y)$ is an open set in an affine $(m-2n)$--dimensional subspace of $V$, which is a translate of the vertical space $\ker(\mathrm{d}\varphi_x)$ for any $x\in \varphi^{-1}(y)$. The quantity $\sigma(y)\in V/\ker(\mathrm{d}\varphi_x)$ is the {\em position vector of $\varphi^{-1}(y)$} and does not depend on the choice of $x$.

\begin{thm}
\label{thm:foli-holo}
Notation as above. The vector bundle $Q^\varphi$ is a holomorphic isotropic quotient and the section $\sigma$ is holomorphic.
\end{thm}

\proof
The holomorphy of $Q^\varphi$ is equivalent, via the universal property, wth the holomorphy of the map $\omega$.
Since the differential of the map $\varphi$ satisfies $\mathrm{d}\varphi\circ F=J_N\circ\mathrm{d}\varphi$, it suffices to prove that $\Omega$ is holomorphic with respect to the $f$--structure $F$ i.e.  $\mathrm{d}\Omega:TV\to \Omega^{-1}T\mathbb F$ satisfies $\mathrm{d}\Omega\circ F=J_{\mathbb F}\circ \mathrm{d}\Omega$.

Let us assume for the moment that we have the following property (which will be proved a little later):
\begin{equation}
\label{eqn:N}
N_F(X,Z)=0\mbox{ for any horizontal }X\mbox{ and any vertical }Z.
\end{equation}
 
%%%%%%%%%%%%%%%%%%%%
%
%For the holomorphy of $\omega$ we only need to prove it along a section $s$, which we will do .
%
%%%%%%%%%%%

From the description of the tangent space at a point of $\mathbb F$, section \ref{sec:homog}, we obtain a vector bundle inclusion $\Omega^{-1}T\mathbb F\subset\mathrm{End}(T_{\mathbb R}U)$. From \cite[Proposition 5.1]{Raw} and \cite[Proposition 5.1]{Raw} it follows that
\begin{equation}
\label{eqn:dOmega}
[\mathrm{d}\Omega(X),F]=\nabla_XF
\end{equation}
for any vector field $X$ on $U$. 
Consider the  splitting of the complexified tanget bundle $T_{\mathbb C}U=T^+U\oplus T^0U\oplus T^-U$ in subbundles corresponding to the eigenvalues $+i$, $0$ and $-i$ of $F$ respectively. With this notation, $T^0U$ is the complexification of the vertical distribution $\ker(F)$ and $T^+U\oplus T^-U$ is the complexification of the horizontal distribution $\im(F)$. In particular, any section in $T^+U$ can be uniquely written as $X-iFX$, where $X$ is a real horizontal vector field on $U$.

The holomorphy of $\mathrm{d}\Omega$ reduces to the compatibility of the corresponding eigenspaces, i.e. it suffices to prove that its complexification maps sections in $T^+U$ to vector fields of type $(1,0)$ on $\mathbb F$ and maps $T^0U$ to $0$.  Lemma \ref{lem:vertical} implies indeed that $\mathrm{d}\Omega(T^0U)=0$.

Note that $\mathrm{ad}(F)$ is injective on $\Omega^{-1}T^+\mathbb F$, and $\mathrm{ad}(F)$ maps $\Omega^{-1}T^+\mathbb F$ to $\Omega^{-1}T^+\mathbb F$, hence, from the description of the holomorphic tangent space to $\mathbb F$ it suffices to verify the three conditions below.

\begin{itemize}
\item 
We verify that $(\nabla_uF)(v)=0$ for any $u\in \Gamma(T^+U)$ and $v\in \Gamma(T^+U)$.

Write $u=X-iFX$ and $v=Y-iFY$. After a small computation using the sPHH assumption, the condition is equivalent to $\nabla_XFY+\nabla_{FX}Y$ being horizontal. We evaluate 
$
g_M(\nabla_XFY+\nabla_{FX}Y,Z)
$ 
for an arbitrary vertical vector field $Z$ and, using the skew-symmetry of $F$, we find that it vanishes for any $Y$ if and only if
\[
F\nabla_XZ=\mathscr{H}\nabla_{FX}Z
\]
for any $Z$.
From the sPHH assumption, the latter equality is equivalent to $F[X,Z]=\mathscr{H}[FX,Z]$, or $F^2[X,Z]=F[FX,Z]$, which follows (\ref{eqn:N}). 

\item
We verify that $(\nabla_uF)(Y)\in \Gamma(T^+U)$ for any $u\in \Gamma(T^+U)$ and $Y$ real, vertical vector field.

Write $u=X-iFX$. Then $(\nabla_uF)(Y)\in \Gamma(T^+U)$  if and only if $F\nabla_XY=\mathscr{H}\nabla_{FX}Y$. As we have already seen, this is implied by $N(X,Y)=0$. 

\item 
We verify that $(\nabla_uF)(v)\in \Gamma(T^0U)$ for any $u\in \Gamma(T^+U)$ and $v\in \Gamma(T^-U)$.

Write $u=X-iFX$ and $v=Y+iFY$. From the sPHH condition, we compute $(\nabla_uF)(v)$ and find it equal to
\[
\mathscr{V}\nabla_XFY-\mathscr{V}\nabla_{FX}Y-i\mathscr{V}\nabla_{FX}FY-i\mathscr{V}\nabla_XY\in \Gamma(T^0U).
\]
\end{itemize}

We prove the remaining relation (\ref{eqn:N}) together with the holomorphy of $\sigma$, following closely \cite[Lemma 2.6]{BW}. The statement is local and hence we may assume that $Q^\varphi$ is trivial. Then $\sigma$ can be interpreted as a map $\sigma:N\to V$ such that $\sigma(y)$ is horizontal for any $y$. 

 It suffices to prove (\ref{eqn:N}) along any section of $\varphi$. Any (local) section $s:N\to U$ of $\varphi$ can be decomposed as $s(y)=\sigma(y)+Z(y)$ with $Z(y)\in\varphi^{-1}(y)$; obviously $Z$ depends on $s$, whereas $\sigma$ does not. If $Y$ is a vector field on $N$ then, since $\varphi$ is horizontally holomorphic, we obtain $\mathscr{H}\mathrm{d}s(J_NY)=F\mathrm{d}s(Y)$ and hence 
\[
\mathscr{H}\mathrm{d}\sigma(J_NY)+\mathscr{H}\mathrm{d}Z(J_NY)=
F\{\mathrm{d}\sigma(Y)+\mathrm{d}Z(Y)\}.
\]
As in the case $n=1$ \cite{BW}, $\mathscr{H}\mathrm{d}Z(Y)$ identifies  with $(\mathrm{d}\Phi\circ\mathrm{d}\omega)(Y)(Z)$. Using relation (\ref{eqn:dPhi}) in the equality above, we obtain 
\[
\mathscr{H}\mathrm{d}\sigma(J_NY)+(\mathrm{d}\Phi\circ\mathrm{d}\omega)(J_NY)(Z)=
F\mathrm{d}\sigma(Y)+F(\mathrm{d}\Phi\circ\mathrm{d}\omega)(Y)(Z).
\]

For $\varepsilon$ small, $s_\varepsilon=\sigma+\varepsilon\cdot Z$ is also a section, and applying the same calculations for $s_\varepsilon$ and identifying the coefficients we find that 
\begin{equation}
\label{eqn:sigma}
\mathscr{H}\mathrm{d}\sigma(J_NY)=F\mathrm{d}\sigma(Y)
\end{equation}
and
\begin{equation}
\label{eqn:Z}
(\mathrm{d}\Phi\circ\mathrm{d}\omega)(J_NY)(Z)=F(\mathrm{d}\Phi\circ\mathrm{d}\omega)(Y)(Z).
\end{equation}

Since $Z$ is vertical we have $(\mathrm{d}\Phi\circ\mathrm{d}\omega)(Y)(Z)=\mathrm{d}\omega(Y)(Z)$ and $(\mathrm{d}\Phi\circ\mathrm{d}\omega)(J_NY)(Z)=\mathrm{d}\omega(J_NY)(Z)$, section \ref{sec:homog} and hence, using the horizontal holomorphy of $s$, we obtain the following equality
\[
\mathrm{d}\Omega(FX)(Z)=F\mathrm{d}\Omega(X)(Z)\mbox{ for }X\mbox{ along }s(N)
\]
which implies $F\mathrm{d}\Omega(FX)(Z)=F^2\mathrm{d}\Omega(X)(Z)$ i.e. $F[\mathrm{d}\Omega(X),F](Z)=[\mathrm{d}\Omega(FX),F](Z)$, as $Z$ is vertical. This shows, after applying (\ref{eqn:dOmega}) that $F(\nabla_XF)(Z)=(\nabla_{FX}F)(Z)$ i.e. $F^2\nabla_XZ=F\nabla_{FX}Z$, which is equivalent, from the sPHH assumption, to $N_F(X,Z)=0$. We have proved relation (\ref{eqn:N}) along $s(N)$.

The relation (\ref{eqn:sigma}) proves that $\sigma$ is holomorphic.
\endproof

Summing up, we have shown how to associate to a totally geodesic foliation with extra-assumptions on an Euclidean domain $U$ a holomorphic map from the leaf space and a holomorphic section in the pullback of the universal bundle from the orthogonal Grassmannian. Foliations whose leaves pass all through the origin correspond to the null--section~$\sigma$, case $(a)$ in \cite{BW}, Theorem 2.8, whereas the other foliations correspond to non--trivial sections, case~$(b)$ in loc.cit.

\end{document}